\documentclass[10pt, reqno]{amsart}

\usepackage[T1]{fontenc}
\usepackage{lmodern}
\usepackage[utf8]{inputenc}   % To use accentuated letters \UTF{008E}, \UTF{008F}... directly in the Latex file
\usepackage{hyperref}
\usepackage{amssymb, mathrsfs}
\usepackage{amsmath, amsthm}
\usepackage{enumerate}
\usepackage{xcolor}
\usepackage[a4paper]{geometry}
\usepackage{titlesec}
\usepackage{fancyhdr}
\usepackage{tikz}

\hypersetup{
    colorlinks   = true,
    linkcolor    = blue!55!black,
    citecolor    = blue!55!black,
    urlcolor     = blue!55!black,
    pdftitle     = {Unique continuation properties for the continuous Anderson operator in dimension 2},
    pdfauthor    = {N. Moench},
    pdfsubject   = {Unique continuation, Anderson operator, quasiconformal mappings},
}

\usepackage{setspace}
\renewcommand{\baselinestretch}{0.99}

\numberwithin{subsection}{section}
\numberwithin{subsubsection}{subsection}
\numberwithin{equation}{section} % Une equation p.q. est la q-eme equation de la section p.

\renewcommand\thesection       {\arabic{section}}
\renewcommand\thesubsection    {\thesection{\boldmath $.$}\arabic{subsection}}
\renewcommand\thesubsubsection    {\thesection{\boldmath $.$}\arabic{subsection}{\boldmath $.$}\arabic{subsubsection}} %{{\it \thesection}{\boldmath $.$}{\it \arabic{subsection}}{\boldmath $.$}{\it \arabic{subsubsection}}}

\titleformat{\section}[block] %\titleformat{\section}[block]
{\filcenter\normalfont\sffamily\bfseries\Large}  % {\filcenter\normalfont\sffamily\bfseries\Large}
{{\hspace{0cm}}\thesection \hspace{0.2em} --\vspace{0cm}}{0.5em}{} % {{\hspace{-0.87cm}}\thesection \hspace{0.2em} --\vspace{0cm}}{0.5em}{}

\titleformat{\subsection}[runin]
{\filcenter\normalfont\sffamily\bfseries\large}  % {\filright\normalfont\sffamily\bfseries\large}
{{\hspace{0cm}}\thesubsection \hspace{0.15em} -- \vspace{0.1cm}}{.2em}{}   %{\thesubsection \hspace{0.2em} -- \vspace{0.1cm}}{.2em}{}   %{\hspace{-0.7cm}\thesubsection \hspace{0.5em} \vspace{0.3cm}}{.5em}{}  
\titlespacing{\subsection}{-0pc}{1.5ex plus .1ex minus .2ex}{0pc}   %\titlespacing{\subsection}{-0pc}{1.5ex plus .1ex minus .2ex}{0pc}

\titleformat{\subsubsection}[runin]
{\filcenter\normalfont\sffamily\bfseries}   %{\normalfont\sffamily\bfseries}
{\filright\sffamily{\hspace{0cm}}\thesubsubsection\hspace{0.2em} --}{.5em}{}\titlespacing{\subsection}{-0pc}{1.5ex plus .1ex minus .2ex}{0pc}

\newtheoremstyle{mystyle}
{3pt}               %space above
{3pt}               %space below
{\it }                      %bodyfont
{}                      %indent
{\bfseries}      % {\sffamily\bfseries}       %headfont
{}                      %punctuation
{0.5em}                 %space after head
{\hspace{0cm}\textit{#2 --} {\hspace{-0.02cm}}\textit{#1}}   %   {#1 #2{\hspace{0.2cm}--\hspace{-0.2cm}}  }   %{\llap{#2 }#1{\hspace{0.2cm}--}}

\theoremstyle{mystyle}

\newtheorem{thm}{Theorem.}   % \newtheorem{thm}{Theorem \hspace{-0.13cm}{\boldmath $.$}}
\newtheorem*{thm*}{Theorem.}

\newtheorem{lem}[thm]{{Lemma}. }%[section]
\newtheorem{prop}[thm]{{Proposition.}}%[chapter]
\newtheorem*{rem*}{Remark.}

\newtheoremstyle{mystyle3}
{3pt}               %space above
{3pt}               %space below
{\it }                      %bodyfont
{}                      %indent
{\bfseries}      % {\sffamily\bfseries}       %headfont
{}                      %punctuation
{0.5em}                 %space after head
{\hspace{-0.8cm}{\textbf{\textit{#2}} --} {\hspace{-0.02cm}}{\textbf{\textit{#1}}}}

\theoremstyle{mystyle3}

\newtheoremstyle{mystyle2}
{3pt}               %space above
{3pt}               %space below
{\it }                      %bodyfont
{}                      %indent
{\bfseries}    % {\sffamily\bfseries}             %headfont
{}                      %punctuation
{0.5em}                 %space after head
{\llap{#2 }{\textbf{\textit{#1}}{\hspace{0.2cm}--}}}

\theoremstyle{mystyle2}

\newtheorem*{definition*}{Definition}
\newtheorem*{theorem*}{Theorem}
\newtheorem*{Remark*}{Remark}
\newtheorem*{lem*} {Lemma}
\newtheorem*{defn*} {Definition}
\newtheorem*{prop*} {Proposition}
\newtheorem*{cor*} {Corollary}

\newcommand{\ssk}{\smallskip}

\renewcommand{\epsilon}{\varepsilon}
\newcommand{\eps}{\epsilon}

\newcommand{\mcD}{\mathcal{D}}

\newcommand{\mcH}{\mathcal{H}}

\newcommand{\mcN}{\mathcal{N}}
\newcommand\mcO{\mathcal{O}}

\newcommand{\bfC}{\mathbf{C}}

\newcommand{\bfN}{\mathbf{N}}

\newcommand{\bfR}{\mathbf{R}}

\newcommand{\bfT}{\mathbf{T}}

\makeatletter
\newcommand*{\defeq}{\mathrel{\rlap{%
                     \raisebox{0.3ex}{$\m@th\cdot$}}%
                     \raisebox{-0.3ex}{$\m@th\cdot$}}%
                     =}
\makeatother

\makeatletter
\newcommand*{\eqdef}{=\mathrel{\rlap{%
                     \raisebox{0.3ex}{$\m@th\cdot$}}%
                     \raisebox{-0.3ex}{$\m@th\cdot$}}%
                     }
\makeatother

\newcommand{\norme}[1]{\left\lVert #1 \right\lVert}

\newcommand{\dd}{\mathrm{d}}
\renewcommand{\phi}{\varphi}
\renewcommand{\div}{\operatorname{div}}

\begin{document}

\begin{center}
{\LARGE\sffamily{\textbf{Unique continuation and nodal sets for continuous Anderson operators}   \vspace{0.5cm}}}
\end{center}

\begin{center}
{\sf  N. MOENCH}
\end{center}

\vspace{1cm}

\begin{center}
\begin{minipage}{0.8\textwidth}
\renewcommand\baselinestretch{0.7} \scriptsize \textbf{\textsf{\noindent Abstract.}} We consider singular continuous Anderson operators $\mcH=-\Delta+\xi$ on the one- and two-dimensional tori, and prove the strong unique continuation property for their eigenfunctions.
The proof uses a ground-state transform to rewrite the eigenvalue equation in divergence form with a Hölder-continuous coefficient, then combines a Carleman estimate in dimension one with a Beltrami equation and Ahlfors–Bers factorization in dimension two. This planar approach further shows that the nodal set of an eigenfunction is locally the quasiconformal image of the zero set of a harmonic function, yielding a Courant-type nodal domain theorem. In dimension one, we upgrade unique continuation into a quantitative spectral inequality which, through the Lebeau–Robbiano method, gives exact null-controllability of the associated parabolic equation from any open subset of the torus.

%\medskip\noindent\textbf{\textsf{Keywords.}} Anderson Hamiltonian, singular SPDEs, unique continuation, quasiconformal mappings, nodal sets, spectral inequality, null-controllability.

%\medskip\noindent\textbf{\textsf{2020 Mathematics Subject Classification.}} 35R60, 35B60, 30C62, 35P05, 93B05.
\end{minipage}
\end{center}

\vspace{0.4cm}
\ssk

%------------------------------------------------------------------------ù
%------------------------------------------------------------------------ù
\section{ Introduction}\label{section_intro}
%------------------------------------------------------------------------ù
%------------------------------------------------------------------------ù

Continuous Anderson operators are random Schrödinger operators of the form
$$
\mcH=-\Delta+\xi,
$$
where the potential $\xi$ is a random spatial noise with insufficient regularity for the product $\xi u$ to be interpreted in the classical sense. In this singular regime, the operator is constructed through a renormalization procedure, using modern theories of singular stochastic partial differential equations such as regularity structures and paracontrolled distributions. These methods provide a rigorous meaning to the Anderson Hamiltonian and allow one to study its spectral properties even when the potential is distribution-valued. In particular, continuous Anderson operators on compact manifolds are known to be self-adjoint, bounded from below, to have compact resolvent, and to satisfy a Weyl law, see \cite{AC,BDM,GUZ,Mouzard}. 

In this paper, we investigate local properties of eigenfunctions of singular Anderson operators. Our main focus is the unique continuation property and its consequences for the geometry of nodal sets. We consider the one and two-dimensional cases and show that, despite the singularity of the underlying potential, eigenfunctions retain a strong form of unique continuation. In dimension two, the planar structure allows us to go further and obtain a quasiconformal description of the nodal set. We also derive a quantitative unique continuation estimate in dimension one and use it to obtain null-controllability for the associated parabolic equation.

\subsection{ Main results.\hspace{0.15cm} }
Throughout the paper, we work with enhanced noises
$$
\Xi=(\xi,\xi^{(2)})\in\mathcal N^\alpha(\mathbf T^d) = C^{\alpha-2}(\bfT^d)\times C^{2\alpha-2}(\bfT^d),
\qquad
\alpha\in(2/3,1),
$$
for $d\in\{1,2\}$, and with the corresponding singular Anderson operator $\mcH$. The construction of the operator and the properties of the enhanced noise are recalled in Section \ref{section_anderson}. 

Our first main result is a strong unique continuation principle.

\begin{thm}\label{thm_cont}
Let $\alpha\in(2/3,1)$, let $\Xi\in\mathcal N^\alpha(\mathbf T^d), d\in\{1,2\}$, and let $\mcH$ be the corresponding singular Anderson operator. Then the eigenfunctions of $\mcH$ satisfy the strong unique continuation property.
\end{thm}
More precisely, if $u$ is an eigenfunction and if, for some $x_0\in\mathbf T^d$,
\begin{equation}\label{eq_vanishinfinity}
\int_{B(x_0,r)}|u|^2=O(r^N)
\end{equation}
for every $N\in\mathbf N$, then $u\equiv0$.

In dimension one, the proof is based on a reduction to a second-order differential equation and on a Carleman estimate. In dimension two, the same ground-state transform leads instead to a uniformly elliptic equation in divergence form. After division by a suitable positive solution of the adjoint equation and introduction of a stream function, the resulting pair can be interpreted as a quasiregular mapping. The Ahlfors-Bers factorization theorem then reduces the problem to the corresponding property for holomorphic functions.

The two-dimensional argument yields additional information on the nodal set.

\begin{thm}\label{thm_quasiconformal}
Let $\Xi\in\mathcal N^\alpha(\mathbf T^2)$, and let $u$ be an eigenfunction of the corresponding singular Anderson operator. Around each zero $x_0$, the nodal set of $u$ is locally the quasiconformal preimage of the zero set of a non-trivial harmonic function.
\end{thm}
In particular, the local geometry of the nodal set is inherited from the nodal set of a harmonic function. If the associated holomorphic function has a zero of order $m$, then locally its real part has $2m$ nodal arcs meeting at the zero. Thus the quasiconformal representation provides a precise description of the possible local singularities of the nodal set.

As a consequence of strong unique continuation, we also obtain a Courant-type nodal domain theorem. With the eigenvalues indexed increasingly as
$$
\lambda_0\leq\lambda_1\leq\cdots,
$$
an eigenfunction $u_n$ associated with $\lambda_n$ has at most $n+1$ nodal domains. We refer to \cite{Alles_courant} for the argument proving the Courant theorem from the strong unique continuation property.

For every non-empty open subset $\omega\subset\mathbf T$, we prove the spectral inequality
$$
\lVert P_\lambda u\lVert_{L^\infty(\mathbf T)}
\leq
e^{C\sqrt{\lambda-\lambda_0}}
\lVert P_\lambda u\lVert_{L^\infty(\omega)},
$$
where $P_\lambda$ is the orthogonal projector onto the subspace $E_\lambda = \text{Span}\{u_k, \, \lambda_k\leq\lambda   \}$. The proof relies on an interpolation inequality for solutions of the divergence form equation obtained after the ground-state transform, together with the three-circles theorem and quasiconformal distortion estimates.
Finally, we apply this estimate to the parabolic equation
$$
(\partial_t+\mcH)g=f1_\omega
\quad\text{on }\mathbf T\times(0,T),
\qquad
g(0,\cdot)=g_0.
$$
Using the spectral decomposition of $\mcH$ and the Lebeau Robbiano method, we obtain exact null-controllability from any non-empty open subset of the torus.

\begin{thm}\label{thm_control}
Let $\Xi\in\mathcal N^\alpha(\mathbf T)$, and let $H$ be the corresponding Anderson operator. For every non-empty open subset $\omega\subset\mathbf T$ and every $T>0$, the equation
$$
(\partial_t+\mcH)g=f1_\omega
$$
is exactly null-controllable in time T. More precisely,
for every
$
g_0\in L^2(\mathbf T),
$
there exists a control
$
f\in L^2((0,T)\times\omega)
$
such that the corresponding solution satisfies
$
g(T,\cdot)=0.
$
\end{thm}
Thus, the quantitative unique continuation estimate provides a direct bridge between the spectral properties of the singular Anderson operator and the control of the associated parabolic dynamics.

\subsection{ Main ideas.\hspace{0.15cm}  }

The proofs are based on a simple structural observation. Let $u_0$ be the ground state of $H$, which is known to be positive, and write
$$
u_0=e^Z.
$$
After conjugating by $u_0$, the singular potential is removed from the equation and replaced by the uniformly elliptic coefficient \(e^{2Z}\). Any eigenfunction becomes a solution of a divergence-form equation with a positive Hölder continuous coefficient. In this way, the singular nature of the original Schrödinger operator is replaced by a rough but bounded elliptic coefficient.

In dimension one, this reduction allows us to construct a suitable primitive of the derivative of the conjugated eigenfunction. The resulting function satisfies a second-order equation to which a Carleman estimate can be applied. A weighted Caccioppoli estimate first shows that infinite order vanishing is preserved under the reduction. The Carleman argument then implies that the function must vanish identically.

The situation is different in dimension two. The divergence-form equation admits a local stream function $s$, and the corresponding complex-valued function $w=v+is$ satisfies a Beltrami equation
$$
\bar\partial w=\mu\partial w
$$
with $\lVert \mu\lVert_{L^\infty}<1$. The Ahlfors-Bers factorization theorem provides a decomposition
$$
w=h\circ\chi,
$$
where $h$ is holomorphic and $\chi$ is quasiconformal. Strong unique continuation therefore follows from the corresponding property for holomorphic functions, while the same factorization gives the local structure of the nodal set.

The quantitative argument in dimension one relies on the same divergence-form structure, after lifting the equation to the cylinder $\mathbf T\times\mathbf R$. Quasiconformal distortion allows us to compare ordinary balls with the deformed balls arising from the conjugating map, while Hadamard's three-circles theorem provides the required interpolation estimate. This ultimately yields the spectral inequality used in the controllability argument.

\subsection{ Relation with previous work.\hspace{0.15cm}  }

Unique continuation for Schrödinger operators with sufficiently regular potentials is classical and can be obtained by Carleman estimates, see for example \cite{Aronzajn,Jerison,LebeauRousseau}. In the singular setting considered here, the potential is distribution-valued and the corresponding arguments cannot be applied directly. The ground-state transform provides a way around this difficulty by replacing the singular potential with a rough elliptic coefficient.

The use of quasiconformal mappings to study unique continuation for planar divergence-form equations goes back to the classical theory of generalized analytic functions, see \cite{Astala}\cite{Vekua}\cite{Schultz}. The present work combines this planar structure with the ground-state transform for singular Anderson operators.

For the control problem, the connection between quantitative unique continuation, spectral inequalities, and null-controllability originates in the Lebeau-Robbiano method \cite{LebeauRobianno}. We follow the approach of \cite{Alles_control}, where quasiconformal techniques are used to establish spectral inequalities for one-dimensional divergence-form operators.

The present work is restricted to the range
$$
\alpha\in(2/3,1),
$$
which contains the regimes in which the required constructions of the Anderson operator are currently available in the form used here. It is natural to expect that the arguments should extend, after suitable modifications of the renormalized structure, to the full subcritical range $\alpha\in(0,1)$. Extending the unique continuation theory to dimension three appears to require a different approach, since the quasiconformal factorization used in dimension two is intrinsically planar.

\subsection{ Organization of the paper.\hspace{0.15cm}  }

Section \ref{section_anderson} recalls the basic properties of the Anderson operator needed in the paper and introduces the ground-state transform. In particular, we rewrite the eigenvalue equation in divergence form and collect the regularity properties of the resulting coefficient.

Section \ref{sectionintermediate} is devoted to strong unique continuation in dimension one. After reducing the eigenvalue equation to a suitable second-order equation, we establish the preservation of infinite-order vanishing and apply a Carleman estimate.

Section \ref{section2} treats the two-dimensional case. We introduce a positive solution of the adjoint equation, construct a stream function, and derive a Beltrami equation. The Ahlfors-Bers factorization then yields both strong unique continuation and the local quasiconformal description of the nodal set.

In Section \ref{section_control}, we develop the quantitative unique continuation estimate on the one-dimensional torus and prove the spectral inequality. Finally, we apply the Lebeau-Robbiano method to obtain exact null-controllability for the parabolic equation associated with the Anderson operator.

%---------------%------------------------------------------------------------%
\section{The Anderson operator and the ground-state transform} \label{section_anderson}
%---------------------------------------------------------------------------

We first recall the properties of the singular Anderson operator that will be used throughout the paper. The main point of this section is to introduce the ground-state transform, which allows us to rewrite the eigenvalue equation in divergence form with a Hölder continuous coefficient. This reduction will be the common starting point for the arguments in dimensions one and two.

\subsection{ Singular Anderson operators.\hspace{0.15cm}  }

Let $d\in\{1,2\}$ and $\alpha\in(2/3,1)$. We consider an enhanced noise
$\Xi=(\xi,\xi^{(2)})\in\mathcal N^\alpha(\mathbf T^d)$, where $\mathcal N^\alpha(\mathbf T^d)$ is the closure of the set of smooth enhanced data
$$\big\{(\xi,\xi\Delta^{-1}\xi-c), \quad \xi\in C^\infty(\mathbf T^d),\, c\in\mathbf R\big\}$$
in $C^{\alpha-2}(\mathbf T^d)\times C^{2\alpha-2}(\mathbf T^d)$. The second component $\xi^{(2)}$ encodes the singular product $\xi\Delta^{-1}\xi$, which is not classically well-defined in the present regularity regime.

The corresponding continuous Anderson operator is obtained by approximation. Given smooth functions $\xi_n$ such that
$(\xi_n,\xi_n\Delta^{-1}\xi_n-c_n)\to\Xi$,
we consider the classical operators
$\mcH_n=\Delta+\xi_n-c_n$. The constants $c_n$ are the renormalization constants: they compensate for the diverging deterministic part of the singular product $\xi_n\Delta^{-1}\xi_n$ as the regularization is removed. In other words, renormalization consists in subtracting the appropriate divergent counterterm before taking the limit. The operators $\mcH_n$ then converge, in the norm resolvent sense, to a self-adjoint operator $\mcH=\mcH_\Xi$, which is the continuous Anderson operator associated with $\Xi$, see \cite{AC,BDM,GIP}.

The operator $\mcH$ is bounded from below and has compact resolvent. Its spectrum is therefore discrete, and we write
$\lambda_0\leq\lambda_1\leq\cdots$, with $\lambda_n\to+\infty$.
Moreover, the ground state $u_0$ associated with $\lambda_0$ can be chosen strictly positive, see \cite[Corollary 16]{BDM}. This positivity will be the key ingredient in the ground-state transform developed below.

\subsection{The ground state transform. \hspace{0.15cm}  }

Let $u_0$ be the positive ground state introduced above. Since $u_0$ is positive and belongs to $C^\alpha(\mathbf T^d)$, we may write
$$
u_0=e^Z
$$
for a function
$
Z\in C^\alpha(\mathbf T^d).
$
Since $\mathbf T^d$ is compact, the coefficient
$
a(x):=e^{2Z(x)}
$
is uniformly positive and belongs to $C^\alpha(\mathbf T^d).$

The conjugated operator is defined by $$\widetilde \mcH v=u_0^{-1}\mcH(u_0v)$$ on the domain $D(\widetilde H)=u_0^{-1}D(H)$. Since $u_0$ is continuous and strictly positive on the compact torus, multiplication by $u_0$ is bounded and invertible on $L^2(\mathbf T^d)$. For smooth approximations, a direct computation using $\mcH u_0=\lambda_0u_0$ gives
\begin{equation}\label{eq_conjugez}
\widetilde{\mcH}w= -e^{-2Z}\div\left(e^{2Z}\nabla w\right)
+\lambda_0w.
\end{equation}
Passing to the limit by the resolvent convergence of the approximating operators yields the same identity in the weak sense for the singular operator. The eigenvalue problem rewrites as 
\begin{equation}\label{eq_eigenvalueproblem}
\div(e^{2Z}\nabla w) + e^{2Z}(\lambda-\lambda_0)w = 0,
\end{equation}
in the weak sense.

This reformulation is the main structural ingredient of the paper. The singular potential $\xi$ no longer appears explicitly in the equation. Instead, it is encoded in the coefficient
$
a=e^{2Z},
$
which is positive, bounded, and Hölder continuous.

The advantage of \ref{eq_eigenvalueproblem} is that it places the eigenvalue problem within the theory of elliptic equations in divergence form. In dimension one, we will use this representation to reduce the problem to a second-order equation and then apply a Carleman estimate. In dimension two, it allows us to exploit the planar structure of uniformly elliptic divergence form equations and to construct a stream function. This leads to a Beltrami equation and, ultimately, to a quasiregular representation.

We will repeatedly use the fact that multiplication by $u_0$ and by $u_0^{-1}$ preserves local vanishing properties, since both functions are continuous and bounded from above and below by positive constants. In particular, the infinite-order vanishing condition
$$
\int_{B(x_0,r)}|u|^2\dd x=O(r^N)
\qquad\text{for every }N\in\mathbf N
$$
is equivalent, locally, to the corresponding condition for the conjugated function $v=u/u_0$.

%..................................%
\section{Strong unique continuation in dimension 1} \label{sectionintermediate}
%...................................%

In this section we prove Theorem \ref{thm_cont} in dimension 1 by applying the Carleman method after an adequate change of variables. We set an extended data $\Xi\in {\boldsymbol{\mcN}}^\alpha(\bfT)$ and write $\mcH$ for the corresponding Anderson operator. We also assume we are given $\tilde{u}$, an eigenfunction of $\mcH$ that vanishes at infinite order at some point $x_0\in\bfT$ in the sense of \eqref{eq_vanishinfinity}.
We let $u=\tilde{u}/u_0$

As $u$ is of class $C^1$, we define for $x$ close to $0$ the function $v$ by
$$
v(x) = \int_0^x e^{2Z(s)}u'(s)\dd s.
$$
The function $v$ satisfies the equation
\begin{equation}\label{eq_edov}
-v'' = -(e^{2Z}u')' = e^{2Z}(\lambda-\lambda_0)u, 
\end{equation}
so $v$ is in the Sobolev space $H^2$. Let us show that $v$ vanishes at infinite order around 0. To do so we use the following Caccioppoli type inequality, that is valid in any dimension.

 \begin{prop}\label{prop_caccioppoli}
For any $Z\in C^\alpha$, any $w\in \mcD( \widetilde\mcH)$ and $r>0$, we have the estimate
$$
\int_{B(0,r/2)} e^{2Z} |\nabla w|^2 \lesssim \frac{1}{r^2}\int_{B(0,r)} e^{2Z} |w|^2 +  r^2\int_{B(0,r)} e^{-2Z}|\div (e^{2Z}\nabla w)|^2 .
$$     
\end{prop}
\begin{proof}
The proof follows the same path as for the classical Caccioppoli estimate. Let $\theta$ be a smooth non-negative cut-off function vanishing outside $B(0,1)$ and such that $\theta=1$ on $B(0,1/2)$, and set $\tilde{\theta}(x) \defeq \theta(x/r)$. We have by integration by parts
\begin{align*}
\int_{B(0,r)} \tilde{\theta}^2 e^{2Z} |\nabla w|^2
   &= - \int_{B(0,r)} w \, \div \big(\tilde{\theta}^2 e^{2Z}\nabla w\big)
   \\
   &= -\int_{B(0,r)} \tilde{\theta}^2 w  \div(e^{2Z}\nabla w)  -  2\int_{B(0,r)} \tilde{\theta} e^{2Z} w \nabla w \cdot \nabla\tilde{\theta}
   \\
   &=:I_1+I_2.
\end{align*}
We use the following Young inequality valid for any $a,b\in\bfR$ and $\eta>0$
$$
ab\leq \frac{1}{\eta} a^2+\eta b^2.
$$
 This gives 
$$
I_1\leq \frac{1}{r^2}\int_{B(0,r)}e^{2Z} |w|^2 + r^2\int_{B(0,r)}e^{-2Z}|\div(e^{2Z}\nabla w)|^2,
$$
and 
\begin{align*}
I_2&\leq \frac{1}{\eps r^2}\int_{B(0,r)}e^{2Z} |w|^2 + \eps r^2\int_{B(0,r)}e^{2Z}\tilde{\theta}^2|\nabla\tilde{\theta}|^2|\nabla w|^2
\\
&\leq  \frac{1}{\eps r^2}\int_{B(0,r)}e^{2Z} |w|^2 + \eps r^2 \big\lVert \nabla\tilde{\theta}\big\lVert_{L^\infty}^2\int_{B(0,r)}\tilde{\theta}^2 e^{2Z}|\nabla w|^2
\end{align*}

As $\big\lVert\nabla\tilde{\theta}\big\lVert_{L^\infty}\lesssim 1/r$, choosing $\eps$ small enough, depending only on $\theta$, one can absorb the integral $\eps\norme{\nabla\tilde\theta}_{L^\infty}^2\int_{B(0,r/2)}e^{2Z} \tilde{\theta}^2 |\nabla w|^2$ into the left-hand side. We conclude the proof by writing
\begin{align*}
    \int_{B(0,r/2)} e^{2Z} |\nabla w|^2 &\leq  \int_{B(0,r)} \tilde{\theta}^2 e^{2Z} |\nabla w|^2
    \\
    &\lesssim \frac{1}{r^2}\int_{B(0,r)}e^{2Z} |w|^2 + r^2\int_{B(0,r)} e^{-2Z} | \div (e^{2Z}\nabla w )|^2.
\end{align*}
\end{proof}

\begin{lem}
    The functions $v$ and $v'$ vanish at infinite order at 0.
\end{lem}
\begin{proof}
    From Proposition \ref{prop_caccioppoli} and using the fact that $u$ is an eigenfunction of $\widetilde{\mcH}$, we have
 \begin{align*}
     \int_{B(0,r/2)}  e^{2Z}|u'|^2 &\lesssim \frac{1}{r^2}\int_{B(0,r)} e^{2Z}|u|^2 + r^2\int_{B(0,r)} e^{2Z} | u|^2  
     \\
     &\lesssim \frac{1}{r^2}\int_{B(0,r)} |u|^2 ,
 \end{align*}
As $e^{2Z}$ is bounded from below, this ensures $\int_{B(0,r/2)}  e^{2Z}|\nabla u|^2 = \mcO(r^N)$ for any $N\in\bfN$. Then $u'$, and hence $v'$, vanish at infinite order at $0$. And from the Cauchy-Schwarz inequality,
 $$
 |v(x)|^2 \lesssim   \Big| \int_0^x |u'|^2 \Big|  = \mcO(|x|^N),
 $$
 for any $N\in\bfN$.
\end{proof}

We have $u(x) = \int_0^x  e^{-2Z(s)}v'(s) \dd s$. Then, from Equation \eqref{eq_edov} and the assumption that $u$ is an eigenfunction of the conjugated operator given by \eqref{eq_conjugez}, it follows that
\begin{equation}\label{Eqv}
v''(x) =  - \big(\lambda - \lambda_0\big)e^{2Z(x)}\int_0^x e^{-2Z(s)}v'(s) \dd s.
\end{equation}

\smallskip

To prove the unique continuation property we will use the following Carleman estimate from Aronszajn's work on strong unique continuation for the Laplace operator.
\begin{thm}(Aronszajn \cite{Aronzajn})
There exists a constant $C$ such that for any $r\in(0,1)$, any smooth $w$ with support included in $B(0,r)\backslash \{0\}$ and $\beta>0$ we have the inequality
$$
\int_{|x|<r} \big(|w|^2 + |\nabla w|^2 \big) |x|^{-2\beta}\dd x \leq Cr^2 \int_{|x|<r} |\Delta w|^2 |x|^{-2\beta} \dd x.
$$
\end{thm}

\begin{proof}(\emph{Theorem} \ref{thm_cont} in dimension 1.)
We let $\chi$ be a smooth cut-off function that vanishes outside the ball $B(0,r)$ and that is equal to $1$ in $B(0,r/2)$, and that is increasing in $\bfR_-$ and decreasing in $\bfR_+$. We also let $\psi$ be a smooth function vanishing in the ball $B(0,1/2)$ and equal to $1$ outside the ball $B(0,1)$, and set $\psi_j(x) \defeq \psi(jx)$.

We will apply the Carleman estimate to the function $\chi_j v$, where $\chi_j$ is defined as $\chi_j\defeq \psi_j\chi$. This gives 
$$
\int_{|x|<r} \big(|\chi_jv|^2 + |(\chi_j v)'|^2 \big) |x|^{-2\beta}\dd x \leq Cr^2 \int_{|x|<r} | (\chi_jv)''|^2 |x|^{-2\beta} \dd x.
$$
We develop $ (\chi_j v)'$ and $(\chi_jv)''$ using the Leibniz rule. The functions $\nabla\psi_j$ and $\Delta\psi_j$ are supported in the ball $B(0,2/j)$ and verify the estimates $\norme{\nabla \psi_j}_{L^\infty}\lesssim j$ and $\norme{\Delta \psi_j}_{L^\infty}\lesssim j^2$. Using also the fact that $v$ and $v'$ vanish at infinite order, we have
$$
\int_{|x|<r} |\chi \psi_j' v|^2 |x|^{-2\beta}\lesssim  j\int_{|x|<1/j} |v|^2 |x|^{-2\beta} = o_{j \to \infty}(1),
$$
$$
\int_{|x|<r} |\chi \psi_j'' v|^2 |x|^{-2\beta}\lesssim  j^2\int_{|x|<1/j} |v|^2 |x|^{-2\beta} = o_{j \to \infty}(1),
$$
and the same with $v'$ replacing $v$ or $\chi'$ replacing $\chi$. Then all the terms where $\psi$ is differentiated go to 0 as $j\to+\infty$, so that letting $j\to+\infty$ we get the inequality
\begin{equation}\label{eq_aronszajn_nocutoff}
\int_{|x|<r} \big(|\chi v|^2 + |(\chi v)'|^2 \big) |x|^{-2\beta}\dd x \leq Cr^2 \int_{|x|<r} | (\chi v)''|^2 |x|^{-2\beta} \dd x.
\end{equation}

\smallskip

From $(\chi v)'' = \chi v'' + 2\chi' v' + \chi'' v $, we have
$$
 \int_{|x|<r} |(\chi v)''|^2 |x|^{-2\beta} \dd x \leq 2\int_{|x|<r} \chi^2|v''|^2 |x|^{-2\beta} \dd x + 2\int_{|x|<r} |2\chi'v'+\chi'' v|^2 |x|^{-2\beta} \dd x
$$
Using Equation \eqref{Eqv} and the boundedness of the functions $a$ and $Z$, for some constant $M$
$$
r^2\int_{|x|<r} \chi^2 |v''|^2 |x|^{-2\beta} \dd x \leq Mr^2 \int_{|x|<r} \chi^2 \Big(\int_0^x |v'| \Big)^2 |x|^{-2\beta}
$$
We will use the following weighted Hardy inequality that is valid for any $\beta>1/2$ and measurable function $f$
$$
\int_\bfR  |y|^{-2\beta} \Big(\int_0^y f   \Big)^2 \dd y\leq \frac{1}{(\beta-1/2)^2} \int_\bfR f(x)^2|x|^{-2\beta+2} \dd x.
$$
From the monotonicity of $\chi$, we obtain for $\beta\geq1$
\begin{align*}
r^2\int_{|x|<r} \chi^2 \Big(\int_0^x |v'| \Big)^2 |x|^{-2\beta} \dd x &\leq 
r^2\int_{|x|<r}  \Big(\int_0^x \chi |v'| \Big)^2 |x|^{-2\beta}
\\
&\leq 4 r^2\int_{|x|<r} (\chi v')^2 |x|^{-2\beta+2} \dd x
\leq 4 r^4 \int_{|x|<r} (\chi v')^2 |x|^{-2\beta} \dd x
\end{align*}
 Choosing $r$ small, one can absorb this last term into the left-hand side of the Carleman inequality \eqref{eq_aronszajn_nocutoff} to obtain 
 \begin{align*}
 \int_{|x|<r} \big(|\chi v|^2 + |(\chi v)'|^2 \big) |x|^{-2\beta}\dd x &\leq C' r^2 \int_{|x|<r} |2\chi'v'+\chi'' v|^2 |x|^{-2\beta} \dd x 
 \end{align*}
 As $\chi=1$ on $(-r/2,r/2)$ and as the support of $\chi'$ and $\chi''$ are both included in $\{x, \enskip r/2<|x|<r \}$, we obtain 
 \begin{align*}
 \int_{|x|<r/2} \big(| v|^2 + | v'|^2 \big) |x|^{-2\beta}\dd x &\leq C'r^2 \int_{|x|<r} |2\chi'v'+\chi'' v|^2 |x|^{-2\beta} \dd x 
 \\ 
 &\leq C'' \int_{r/2<|x|<r} \big(|v'|^2+|v|^2 \big) |x|^{-2\beta} \dd x.
 \end{align*}
Then
$$
(r/4)^{-2\beta} \int_{|x|<r/4} \big(| v|^2 + | v'|^2 \big)\dd x\leq C''(r/2)^{-2\beta} \int_{r/2<|x|<r} \big(|v'|^2+|v|^2 \big) \dd x.
$$
It suffices to send $\beta$ to $+\infty$ to obtain $v=0$ in $B(0,r/4)$. We conclude that $v$ is identically zero by a chain of ball argument.
\end{proof}

%------------------------------------------------------%
\section{Strong unique continuation in dimension 2} \label{section2}
%------------------------------------------------------%

%------------------------------------------------------%
\subsection{Quasiregular mappings on the plane. \hspace{0.15cm}} \label{subsect_quasi}
%------------------------------------------------------%

Let $\Omega$ be an open subset of the plane, which we identify with the complex plane $\bfC$. We introduce the Wirtinger derivatives
$$
\partial f = \frac{1}{2}(\partial_x-i\partial_y)f ,\qquad \overline\partial f = \frac{1}{2}(\partial_x+i\partial_y)f.
$$
The class of quasiregular mappings generalizes the class of holomorphic functions and is defined as the maps on the plane with values in $\bfC$ satisfying the Beltrami equation
\begin{equation}\label{eq_beltrami}
    \overline\partial f(z) = \mu(z) \partial f(z),
\end{equation}
for some measurable function $\mu$ with modulus bounded by some constant $\frac{k-1}{k+1}<1$ called the distortion factor. We then say that the map $f$ is $k$-quasiregular. A quasiregular mapping that is homeomorphic is called quasiconformal. These mappings are useful for the study of elliptic equations in divergence form on the plane, see for instance the book \cite{Astala} or the classical treatise \cite{Vekua} on generalized analytic functions. The Ahlfors-Bers representation theorem states that any $k$-quasiregular mapping $f$ on some ball $B$ in the plane factorizes as $f=h\circ\chi$ for a $k$-quasiconformal $\chi$ and a holomorphic function $h$.
Mori's theorem ensures that for any $k$-quasiconformal mapping $\chi$, there exists an exponent $\gamma\in(0,1)$ depending only on $k$ and a constant $C$ such that the following inequality holds for any points $x,y$:
\begin{equation}\label{eq_mori}
\frac{1}{C} |y-x|^{1/\gamma}\leq |\chi(y)-\chi(x)| \leq C|y-x|^\gamma.
\end{equation}

Quasiconformal mappings are useful to study elliptic equations in divergence form on the plane as they relate them to holomorphic functions. We will use them to prove strong unique continuation for Anderson operators in dimension 2, and we will also obtain the following result giving information on the zero set of eigenfunctions.

\subsection{Proof of strong unique continuation property. \hspace{0.15cm} }

We now prove Theorem \ref{thm_quasiconformal}. The idea is to put the equation for the eigenfunction into divergence form and to use tools from the theory of quasiconformal mappings.

We set for the whole section an enhanced noise $\Xi\in{\boldsymbol{\mcN}} ^\alpha(\bfT^2)$ and work with the corresponding Anderson operator. We reproduce the arguments from \cite{Schultz}, which prove strong unique continuation for weak solutions of the divergence-form elliptic equation. We conjugate the Anderson operator $\mcH$ by its ground state $u_0=\exp(Z)$ as was done in Section \ref{sectionintermediate}. We then consider $u$, an eigenfunction of the conjugated operator,
$$
\widetilde{\mcH}w = -e^{-2Z}\div (e^{2Z} \nabla w)+\lambda_0 w.
$$
so that $u$ is a solution of the equation
\begin{equation}\label{eq_udivergenceform}
    \div\big(  e^{2Z} \nabla u\big) + e^{2Z}(\lambda-\lambda_0)u =0.
\end{equation}

We fix a point $x_0\in\bfT^2$ and let $\psi$ be a positive solution in a neighborhood of $x_0$ of the following adjoint equation
$$
L^*\psi= \div (e^{2Z} \nabla\psi)+e^{2Z}(\lambda-\lambda_0)\psi=0,
$$
which exists by the local theory for uniformly elliptic divergence-form equations with bounded measurable coefficients. Then $v\defeq \frac{u}{\psi}$ is in a neighborhood of $x_0$, a weak solution of the divergence-form equation
$$
\div(e^{2Z_1} \nabla v)=0.
$$
with
$$
e^{2Z_1} \defeq e^{2Z}\psi^2.
$$
From the Poincaré lemma, this equation is equivalent to the local existence of a function $s$, called the stream function, such that
\begin{equation}\label{eq_stream}
 e^{2Z_1}\nabla v = *\nabla s \defeq \left(\begin{array}{cc}
     \partial_2 s \\-\partial_1 s,
\end{array}\right)
\end{equation}
with $s(x_0)=0$. We also define $w = v+is$. Note that if $u$ admits a zero of infinite order at $x_0$, then so does $\nabla v$, from the Caccioppoli estimate of Proposition~\ref{prop_caccioppoli}, and then $s$ and $w$ admit a zero of infinite order too. Equation \eqref{eq_stream} rewrites as a Beltrami equation
$$
\overline\partial w = \mu{\partial} {w},
$$
with 
$$
\mu \defeq \frac{e^{2Z_1}-1}{e^{2Z_1}+1} \cdot \frac{\partial_1v+i\partial_2 v}{\partial_1v-i\partial_2 v}
$$
when $\nabla v\ne 0$ and $\mu=0$ elsewhere. As $Z_1$ is a bounded function, there exist a constant $k\geq 1$ such that for all $x\in\bfT^2$ 
$$
|\mu(x)|\leq \frac{k-1}{k+1} <1.
$$
It follows from the Ahlfors-Bers representation theorem that one has, on some ball around $x_0$, the factorization 
\begin{equation}\label{eq_factor}w=h\circ\chi,
\end{equation}for a holomorphic $h$ and a quasiconformal $\chi$. 

Mori's theorem stated in Subsection \ref{subsect_quasi} ensures that if $w$ admits a zero of infinite order, then $h$ admits it as well, and the holomorphic nature of $h$ implies that $h$ is identically zero on the ball where it is defined. From this, Theorem~\ref{thm_cont} in dimension $2$ follows as a consequence. 

We also obtain Theorem \ref{thm_quasiconformal} in the process, because taking the real part of Equation \eqref{eq_factor} gives $v= \operatorname{Re}(h) \circ \chi $, and then locally
$$
\{u=0\} = \chi^{-1}\big( \{\operatorname{Re}(h) = 0\}\big).
$$
The real part of any holomorphic function being a harmonic one, this gives Theorem \ref{thm_quasiconformal}.

%------------------------------------------------------%
\section{A spectral inequality for the Anderson operator in dimension 1. \hspace{0.15cm}} \label{section_control}
%------------------------------------------------------%

Once we have the factorization $w=h\circ\chi$, one can exploit the holomorphic property of $h$ to gain some more quantitative form of unique continuation, which usually takes the form of doubling inequalities. One starts here with Hadamard's three circles theorem, which states that for any holomorphic $h$ around the origin, setting $m(r)=\sup_{|z|\leq r} |h(z)|$, one has the convexity inequality

\begin{equation}\label{eq_hadamard}
    m(r)\leq m(r_1)^\theta m(r_2)^{1-\theta},
\end{equation}
with $r=r_1^\theta r_2^{1-\theta}$ and $\theta\in(0,1)$.

A similar inequality has been proven in \cite{Alles_control} for solutions of the divergence-form elliptic equation on the plane $\div(e^{2Z}\nabla f)=0$ on a disc $B(0,R)$. Suppose we are given such an $f$. With the same arguments as in the last subsection, there exists a stream function $s$ with $s(0)=0$ such that, setting $w=f+is$, one can write $w=h\circ \chi$ for a quasiconformal $\chi$ with $\chi(0)=0$ and $h$ holomorphic. The interpolation inequality is then the following proposition.

\begin{prop} (\cite{Alles_control}) \label{prop_interpolation}
    For $r_1\leq r_2\leq R$ and $r=r_1^\theta r_2^{1-\theta}$ with $\theta\in(0,1)$, there is a constant $C$ such that we have the estimate
    $$
    \sup_{B_\chi(r/2)}|f| \leq  C \sup_{B_\chi(r_1)}|f|^\theta  \sup_{B_\chi(r_2)}|f|^{1-\theta},
    $$
    where $B_\chi(r) = \big\{ z,\enskip |\chi(z)|\leq r \big\}$. 
\end{prop}

\begin{proof}

We redo quickly the proof from \cite{Alles_control}. Decompose the holomorphic function $h$ into real and imaginary parts $h=h_1+ih_2$. We have chosen $s$ such that $s(0)=0$, so $h_2(0)=0$, and the Cauchy-Riemann equations give
$$
h_2(x,y) = \int_0^x \partial_2h_1(t,0)\dd t- \int_0^y\partial_1h_1(x,t)\dd t.
$$

From classical interior estimates on the gradient of harmonic functions, for $r>0$ there is a constant $C$ such that 
$$
\sup_{B(0,r)} |h_2|\leq C \sup_{B(0,2r)} |h_1|.
$$
We can then obtain an equivalence between the size of $h$ and the size of $h_1$:
$$
\sup_{B(0,r)} |h_1| \leq \sup_{B(0,r)} |h| \leq C'\sup_{B(0,2r)} |h_1|.
$$

As $f=h_1\circ \chi$, Hadamard's three circles theorem immediately gives the desired inequality.
\end{proof}

This estimate also applies to solutions of elliptic equations defined on the cylinder $\bfT\times\bfR$ as any function defined on it can be lifted by periodicity to a function defined on $\bfR^2$. We will use this setting to study Anderson operators on the 1 dimensional torus $\bfT$. Another remark is that one can make use of Mori's theorem \eqref{eq_mori} to replace the deformed balls $B_\chi$ with true balls in the inequality by changing their radius. 

For $\theta\in(0,1)$, the following interpolation estimate holds, for some constant $C$: 
\begin{equation} \label{eq_interpolation2}
\norme{f}_{L^\infty(B(0,r/2))} \leq C r^{-\theta/2} \norme{f(0,\cdot)}_{L^2(-r,r)}\norme{f}_{L^\infty(B(0,4r))},
\end{equation}
where $f$ is still a solution of the divergence equation $\div(e^{2Z}\nabla f)=0$. We refer to \cite{Alles_control} for a proof of this result.

\bigbreak

We now prove the spectral inequality \eqref{eq_specineq} for Anderson operators, using the same conjugation as in the previous section and following the method of \cite{Alles_control}. We set here the Anderson operator $\mcH=-\partial_x^2+\xi(x)$ associated to some random noise $\Xi\in{\boldsymbol{\mcN}}_\alpha(\bfT)$. Let $(u_k)_{k\geq 0}$ be an orthonormal basis of $L^2$ consisting of eigenfunctions of $\mcH$ with associated eigenvalues $(\lambda_k)_{k \geq 0}$ sorted in increasing order, and we still write $u_0=e^Z$ here.  We let $P_\lambda$ be
the orthogonal projector onto the subspace $E_{\leq\lambda}=\operatorname{Vect}\big\{ u_k, \enskip \lambda_k\leq\lambda \big\}.$ The spectral inequality of interest takes the following form.

\begin{prop} \label{prop_specineq}
Let $\omega$ be an open subset of $\bfT$. There exists a constant $C$ such that for any $u\in L^2(\bfT) $ one has the inequality
    \begin{equation}\label{eq_specineq}
    \sup_{\bfT} |P_\lambda u| \leq e^{C\sqrt{\lambda-\lambda_0}} \sup_\omega  |P_\lambda u|.
    \end{equation}
\end{prop}

\begin{proof}
Let $\lambda\in\bfR$ and $u\in E_{\leq\lambda}$, written as 
$u(x)=\sum_{\lambda_k\leq \lambda} a_k u_k(x)$. Define on $\bfT\times\bfR$ the function
$$
f(x,y) =\sum_{\lambda_k\leq \lambda} a_k  \cosh \big( \sqrt{\lambda_k-\lambda_0 } \, y \big) \frac{u_k(x)}{u_0(x)}.
$$

The sequence of functions $(u_k/u_0)_{k\geq 0}$ is a basis of eigenfunctions of the conjugated operator $\tilde{\mcH}u= e^{-2Z}\div(e^{2Z}\nabla u)+\lambda_0u$. The function $f$ then satisfies the equation
\begin{equation}\label{eq_pdef1}
\partial_y^2f+ e^{-2Z}\partial_x(e^{2Z}\partial_x f) =0 ,
\end{equation}
and $f(x,0)=\frac{u(x)}{u_0(x)}$.

Setting $\tilde{Z}(x,y)\defeq Z(x)$, Equation \eqref{eq_pdef1} rewrites as 
\begin{equation}\label{eq_pdef2}
\div( e^{2\tilde Z}\nabla f)= 0
\end{equation}
which is an equation in divergence form. Combining Proposition \ref{prop_interpolation} and the interpolation inequality \eqref{eq_interpolation2} gives, for some $\alpha\in(0,1)$ depending on the inradius of $\omega$, the key estimate
\begin{equation}\label{eq_interpolation3}
\norme{f}_{L^\infty(\bfT\times(-1,1))}\lesssim \norme{f}_{L^\infty(\omega\times\{0\} )}^\alpha\norme{f}^{1-\alpha}_{L^\infty(\bfT\times(-2,2))}.
\end{equation}

The function $u_k/u_0$ verifies the equation 
$
-\partial_x(e^{2Z} \partial_x(u_k/u_0)) = e^{2Z}(\lambda_k-\lambda_0)u_k/u_0.
$
Then
$$
\int_{\bfT} e^{2Z} |\partial_x(u_k/u_0)|^2\, \dd x = (\lambda_k-\lambda_0)\int_{\bfT} e^{2Z} |u_k/u_0|^2 \, \dd x,
$$
hence the energy estimate $\norme{\partial_x (u_k/u_0)}^2_{L^2(\bfT)}\lesssim \lambda_k-\lambda_0$. Therefore, for every \(y\in(-2,2)\),
\begin{align*}
\|f(\cdot,y)\|_{H^1(\mathbf T)}
& \lesssim
\sum_{\lambda_k\leq\lambda}
|a_k|\,
\cosh \big(\sqrt{\lambda_k-\lambda_0}\,y\big)
\sqrt{1+\lambda_k-\lambda_0}
\\
&
\lesssim
e^{C\sqrt{\lambda-\lambda_0}}
\Big(\sum_{\lambda_k\leq\lambda}|a_k|^2\Big)^{1/2} = e^{C\sqrt{\lambda-\lambda_0}} \norme{u}_{L^2(\bfT)}.
\end{align*}
Where we used Weyl's law in order to absorb the polynomial factor in $\lambda$ arising from the summation over the eigenmodes into the exponential $e^{C\sqrt{\lambda-\lambda_0}}$.
From the Sobolev injection $H^1(\mathbf T)\hookrightarrow L^\infty(\mathbf T)$, we get
$$
\|f\|_{L^\infty(\mathbf T\times(-2,2))} \lesssim \|f\|_{L^\infty_y((-2,2),H^1_x(\bfT))}
\lesssim 
e^{C\sqrt{\lambda-\lambda_0}}
\norme{u}_{L^2(\bfT)}.
$$
Consequently,
\begin{equation} \label{eq_ineqsobolev}
\|f\|_{L^\infty(\mathbf T\times(-2,2))}
\lesssim
e^{C\sqrt{\lambda-\lambda_0}}
\left\| {u/u_0}\right\|_{L^\infty(\mathbf T)}.
\end{equation}
So that gathering \eqref{eq_interpolation3} with \eqref{eq_ineqsobolev}
$$
\sup_\bfT |u/u_0|\leq  \norme{f}_{L^\infty(\bfT\times(-1,1))}    \lesssim e^{C\sqrt{\lambda-\lambda_0}}\norme{u/u_0}_{L^\infty(\omega)}^\alpha\norme{u/u_0}^{1-\alpha}_{L^2 (\bfT)}
$$
And finally
$$
\norme{u}_{L^\infty(\bfT)}\lesssim e^{C\sqrt{\lambda-\lambda_0}}\norme{u}_{L^\infty(\omega)}.
$$
\end{proof}

The spectral inequality above implies a quantitative observability estimate for finite-dimensional spectral subspaces, which is the main ingredient in the Lebeau-Robbiano method. More precisely, if $E_\lambda$ denotes the spectral subspace associated with eigenvalues smaller than $\lambda$, then Proposition \ref{prop_specineq} gives, for every $v\in E_\lambda$,
$$
\lVert v\lVert_{L^\infty(\mathbf T)}
\leq
e^{C\sqrt{\lambda-\lambda_0}}
\lVert v\lVert_{L^\infty(\omega)}.
$$
Since $v$ belongs to the finite-dimensional spectral subspace $E_\lambda$, we may combine this estimate with the standard local elliptic $L^2$-$L^\infty$ estimates on $E_\lambda$, whose constants grow at most polynomially in $\lambda$. By increasing $C$, this polynomial loss can be absorbed into the exponential factor, and we obtain
$$
\lVert v\lVert_{L^2(\mathbf T)}
\leq
e^{C\sqrt{\lambda-\lambda_0}}
\lVert v\lVert_{L^2(\omega)},
\qquad v\in E_\lambda,
$$
with a constant depending only on the enhanced noise and on the control region $\omega$.

We can then apply the classical Lebeau-Robbiano iteration to the semigroup generated by $H$. Since $H$ is self-adjoint, bounded from below and has compact resolvent, the solution of the homogeneous equation is given by $e^{-tH}$, and the high-frequency part of the solution is exponentially damped by the semigroup. Splitting the solution into low and high frequencies and choosing successively adapted spectral thresholds allows one to construct a control supported in $\omega$ which drives any initial datum $g_0\in L^2(\bfT)$ to zero at time $T>0$. In particular, for every $T>0$, there exists a constant $C_T>0$ such that one can find $f\in L^2((0,T)\times\omega)$ satisfying
$$
\lVert f\lVert_{L^2((0,T)\times\omega)}
\leq C_T\lVert g_0\lVert_{L^2(\mathbf T)}
$$
and such that the corresponding solution of
$$
(\partial_t+\mcH)g=f {\bf 1}_\omega,
\qquad
g(0,\cdot)=g_0,
$$
satisfies $g(T,\cdot)=0$.
This proves the exact null-controllability stated in Theorem \ref{thm_control}.

It is worth emphasizing that the only quantitative input required in the Lebeau-Robbiano argument is the spectral inequality established in Proposition \ref{prop_specineq}. Hence the singular nature of the Anderson operator enters the controllability argument through the ground-state transform and the resulting quantitative unique continuation estimate. Once Proposition \ref{prop_specineq} is available, the remainder of the argument follows the standard spectral decomposition and iteration scheme.

\bigskip 

\noindent \textcolor{gray}{$\bullet$} {\sf N. Moench} -- Univ. Lorraine, CNRS, IECL - BP 70239 F-
54506 Vandœuvre-lès-Nancy, France.   \\
\noindent {\it E-mail}: nicolas.moench@univ-lorraine.fr   

\end{document}